\documentclass[papera4,12pt]{article}
\usepackage[english]{babel}
\usepackage[letterpaper,top=2cm,bottom=2cm,left=3cm,right=3cm,marginparwidth=1.75cm]{geometry}
\usepackage{amsmath}
\usepackage{graphicx}
\usepackage[english]{babel}
\usepackage[colorlinks=true, allcolors=blue]{hyperref }
\usepackage{amsthm}
\usepackage{amsfonts}
\usepackage{hyperref}
\newtheorem{theorem}{Theorem}[section]

\newtheorem{lemma}[theorem]{Lemma}

\numberwithin{equation}{section}
\begin{document}
\begin{center}
{\bf {\Large Two Problems on Narayana Numbers And Repeated Digit Numbers }}
\vspace{8mm}

{ \bf  G. Abou-Elela$^1$, A. Elsonbaty$^2$, and  M. Anwar$^3$}
\vspace{3mm}\\
$^1$ Department of Mathematics, University of Damietta \\
Faculty of science, Egypt \\
$^2$$^3$ Department of Mathematics, University of Ain Shams \\
Faculty of science, Egypt \\
e-mails:\href{mailto:ghadaabouelela@du.edu.eg}{\url{ghadaabouelela@du.edu.eg}}$^1$, \\
\href{mailto:Ahmadelsonbaty@sci.asu.edu.eg}{\url{Ahmadelsonbaty@sci.asu.edu.eg}}$^2$,\\
\href{mailto:mohmmed.anwar@hotmail.com}{\url{mohmmed.anwar@hotmail.com}},\hspace{0.5cm}\href{mailto:mohamedanwar@sci.asu.edu.eg}{\url{mohamedanwar@sci.asu.edu.eg}}$^3$
\end{center}
\vspace{10mm}
\begin{abstract}
This work aims to solve two problems in  Diophantine equation of the Narayana sequence $(OEIS {\color{red}A000930})$. In the first part it’s proved that the Narayana number can not be factored as a product of two repdigit numbers for base $2\leq b\leq50$, except in two cases. In The second it has been  proved that there is a finite number of solutions up to 290 to  express  the  product of two Narayana numbers  as base $b-$ repdigits numbers, $2\leq b\leq50$, the proofs of these results use some number-theoretic technique includes Baker's method of linear forms in logarithms height, and some reduction technique.
\end{abstract}
$\mathbf{key words}:$ Narayana sequence, linear forms in logarithms, and Davenport’s lemma
%\maketitle
\section{Introduction}
Let $\{N_{n}\}_{n\geq0}$ be The Narayana sequence given by $ N_{0}=0,N_{1}=1,N_{2}=1$
and the recurrence relation \\
\begin{equation}
  N_{n}=N_{n-1}+N_{n-3}\quad for\,all\, n\geq 3
  \label{re_rel}
\end{equation}
The first values of $N_{k}$ are $ 0,1,1,1,2,3,4,6,\cdots$. Narayana cow's sequence is a problem similar to the
problem of Fibonacci rabbits as it counts calves produced every four years.
This sequence  $(OEIS {\color{red}A000930})$ appeared for the first time in the book "Ganita Kaumudi" (1365) by
Indian mathematician Narayana Pandita, who gave this sequence his name, and play roles in mathematical developments
such as, finding  the approximate value of the square roots, investigations into the  Diophantine equation $a
x^2+1=y^2$ ({\color{red}Pell's equation}). Narayana  cows sequence, also known as the supergolden sequence and the
real root corresponding to  the solution of the characteristic equation is known as the super golden ratio. In
Pascal's triangle, starting from $ n\geq3$ we find that the sum of its rows with triplicated diagonals is a Narayana
sequence, while the sum of the rows with slops diagonals of 45 degrees express the Fibonacci sequence. This sequence
plays an  important role in cryptography, coding theory, and graph theory.
\\
\\
In this paper\, we determine all the solutions of the Diophantine equation \\
\begin{equation}
 N_{n}N_{m}=[a,\cdots,a]_{b}=a \, (\frac{b^{l}-1}{b-1})
 \label{equ:2}
\end{equation}
in integers $(n,m,l,a,b)$ with $3\leq m\leq n,2\leq b\leq 50,1\leq a\leq b-1 $ and $ l \geq 2$.\\
and the solutions of the Diophantine equation \\
\begin{equation}
 N_{k}=a_1  a_2 \, (\frac{b^{l_1}-1}{b-1}) \, (\frac{b^{l_2}-1}{b-1})
 \label{equ:r}
\end{equation}
in integers $(k,b,a_1,a_2,l_1,l_2)$ with $ 2 \leq l_1\leq l_2$ , $1 \leq a_1\leq a_2 \leq b-1$, $k\geq3$,  and $
b\geq 2$ .\\
\\
Many authors have studied such a diophantine equation, for example, Luca {\color{blue}\cite{Di_An}} showed that $F_{5}=55$ and $L_{5}=11$
are the largest repdigits in the Fibonacci and Lucas sequences respectively, the researchers in {\color{blue}\cite{pfib}} showed that $F_{10}=55$  and  $L_6=18$ it is the largest Fibonacci  and Lucas number respectively that can be expressed as a product of two repdigits, the author in  {\color{blue}\cite{Mahadi}} studied  the sum of three Padovan
numbers  as  repdigits in base 10 and he found them, the researchers in  {\color{blue}\cite{BHOI}} showed that the
only Narayana numbers expressible as sums of two repdigits are $N_{14}=88$  and $ N_{17}=277$.\\
\\
In the following theorem we consider $n\geq3$ because $N_1=N_2=N_3=1$.

\begin{theorem}
    The only  solution to the Diophantine equation {\color{blue}(\ref{equ:r})} are
    \begin{equation*}
      N_{8}=\frac{2^{2}-1}{2-1}\frac{2^2-1}{2-1}={\color{red}[11]_{2}[11]_{2}}
    \end{equation*}
    \begin{equation*}
    and \,  N_{16}=\frac{2^{2}-1}{2-1}\frac{2^6-1}{2-1}={\color{red}[11]_{2}[111111]_{2}}
    \end{equation*}
    \label{th:u}
  \end{theorem}
\begin{theorem}
  Let $ 3\leq m \leq n, b\in\{2,3,..,50\}$, $a \in\{1,..,b-1\}$, and $l\geq2$. If $N_n N_m$is a repdigits in base
  $b $  then the only  solutions are given by
\[  (n,m,l,a,b)\in \left\{
\begin{array}{cccc}
  (5,3,2,1,2),&(6,3,2,1,3),&(9,3,3),&(4,4,2,1,3) \\
 (11,9,2,1,3),&((7,3,2,1,5),&(5,4,5),&(6,4,2,1,7)\\
  (10,5,2,1,7),&(8,3,2,1,8),&(5,5,2,1,8),&(7,4,2,1,11)\\
  (6,5,2,1,11),&(9,3,2,1,12),&(19,6,2,1,13),&(6,6,2,1,15) \\
  (8,4,2,1,17),&(7,5,2,1,17),&(10,3,2,1,18),&(7,6,2,1,23)\\
  (9,4,2,1,25),&(8,5,2,1,26),&(11,3,2,1,27),&(8,6,2,1,35)\\
  (7,7,2,1,35),&(10,4,2,1,37),&(9,5,2,1,38),&(12,3,2,1,40)\\
  (15,10,2,1,49)& & &
\end{array}%
\right \}
\]
  \[
(n,m,l,a,b)\in \left\{
   \begin{array}{ccccc}
     (6,4,2,2,3),&(9,4,3,2,3),&(7,4,2,2,5),&(8,4,2,3,5)\\
     (6,5,2,2,5),&(7,5,2,3,5),&(7,6,2,4,5),&(11,3,2,4,6)\\
      (15,3,3,3,6),&(6,6,2,2,7),&(7,6,2,3,7),&(10,7,3,2,7)\\
       (10,8,3,3,7),&(8,4,2,2,8),&(7,5,2,2,8),&(8,5,2,3,8)\\
       (8,6,2,4,8),&(7,7,2,4,8),&(8,7,2,6,8),&(13,3,2,6,9) \\
     (11,9,3,4,9),&(13,12,4,3,9),&(14,3,2,8,10),&(13,3,2,5,11)\\
     (13,4,2,10,11),&(11,5,2,7,11),&(7,6,2,2,11),&(8,6,2,3,11)\\
     (7,7,2,3,11),&(11,10,3,4,11),&(9,4,2,2,12),&(9,5,2,3,12)\\
     (9,6,2,4,12),&(9,7,2,6,12),&(9,8,2,9,12),&(11,3,2,2,13)\\
     (11,4,2,4,13),&(11,5,2,6,13),&(11,6,2,8,13),&(11,7,2,12,13)\\
     (19,11,4,7,13),&(13,3,2,4,14),&(13,4,2,8,14),&(13,5,2,12,14)\\
     (14,4,2,11,15),&(11,6,2,7,15),&(16,9,3,9,16),&(13,5,2,10,17)\\
     (8,6,2,2,17),&(7,7,2,2,17),&(8,7,2,3,17),&(11,8,2,14,17)\\
     (10,4,2,2,18),&(10,5,2,3,18),&(10,6,2,4,18),&(10,7,2,6,18)\\
     (10,8,2,9,18),&(10,9,2,13,18),&(13,3,2,3,19),&(13,4,2,6,19)\\
     (13,5,2,9,19),&(13,6,2,12,19),&(13,7,2,18,19),&(16,3,2,9,20)\\
     (16,4,2,18,20),&(11,5,2,4,20),&(11,7,2,8,20),&(11,8,2,12,20)\\
     (14,3,2,4,21),&(14,4,2,8,21),&(14,5,2,12,21),&(14,6,2,16,21)\\
     (13,4,2,5,23),&(14,5,2,11,23),&(13,6,2,10,23),&(11,7,2,7,23)\\
     (13,7,2,15,23),&(14,7,2,22,23),&(9,6,2,2,25),&(9,7,2,3,25)\\
     (11,9,2,14,25),&(16,3,2,7,26),&(16,4,2,14,26),&(16,5,2,21,26)\\
     (8,7,2,2,26),&(8,8,2,3,26),&(13,8,2,20,26),&(11,4,2,2,27)\\
     (11,5,2,3,27),&(11,6,2,4,27),&(11,7,2,6,27),&(11,8,2,9,27)\\
     (11,9,2,13,27),&(11,10,2,19,27),&(18,3,2,14,28),&(13,3,2,2,29)\\
     (13,4,2,4,29),&(13,5,2,6,29),&(13,6,2,8,29),&(13,7,2,12,29)\\
     (13,8,2,18,29),&(13,9,2,26,29),&(14,6,2,11,31),&(14,5,2,8,32)\\
     (14,7,2,16,32),&(14,8,2,24,32),&(20,19,4,14,33),&(19,3,2,17,34)\\
     (13,5,2,5,35),&(16,6,2,21,35),&(13,7,2,10,35),&(11,8,2,7,35)\\
     (13,8,2,15,35),&(14,8,2,22,35),&(10,6,2,2,37),&(10,7,2,3,37)\\
     (11,10,2,14,37),&(13,10,2,30,37),&(9,7,2,2,38),&(9,8,2,3,38)\\
     (13,9,2,20,38),&(13,4,2,3,39),&(13,6,2,6,39),&(13,7,2,9,39)\\
     (12,4,2,2,40),&(12,5,2,3,40),&(12,6,2,4,40),&(12,7,2,6,40)\\
     (12,8,2,9,40),&(12,9,2,13,40),&(12,10,2,19,40),&(12,11,2,28,40)\\
     (16,4,2,9,41),&(11,5,2,2,41),&(18,5,2,29,41),&(16,6,2,18,41)\\
     (11,7,2,4,41),&(16,7,2,27,41),&(11,8,2,6,41),&(13,11,2,40,41)\\
     (15,3,2,3,42),&(15,4,2,6,42),&(15,5,2,9,42),&(15,6,2,12,42)\\
     (15,7,2,18,42),&(15,8,2,27,42),&(15,9,2,39,42),&(14,3,2,2,43)\\
   (14,4,2,4,43),&(14,5,2,6,43),&(14,6,2,8,43),&(14,7,2,12,43)\\
     (14,8,2,18,43),&(14,9,2,26,43),&(14,10,2,38,43),&(13,5,2,4,44)\\
     (13,7,2,8,44),&(13,8,2,12,44),&(20,10,3,8,45),&(13,6,2,5,47)\\
     (14,7,2,11,47),&(13,11,2,35,47),&(11,11,2,16,48),&(19,5,2,35,50) %

   \end{array}%
\right \}
\]

\label{th:1}
\end{theorem}

\section{ Preliminary}
\subsection{Narayana sequence}
The characteristic equation corresponding to the third-order linear recurrence relation
${\color{blue}(\ref{re_rel})}$ is  $ x^{3}-x^{2}-1$, this equation  has roots $\alpha$,$\beta$, and
$\gamma=\bar{\beta} $ where \\
\begin{equation*}
  \alpha=\textstyle\frac{2+r_{1}+r_{2}}{6}, \beta=\textstyle\frac{4-(1+\sqrt{-3}) r_{1}-(1-\sqrt{-3})r_{2}}{12}
\end{equation*}
and \\
\begin{equation*}
  r_{1}=\sqrt[3]{116-12\sqrt{93}},r_{2}=\sqrt[3]{116+12\sqrt{93}}
\end{equation*}
Furthermore, the  Bient formula is \\
\begin{equation*}
  N_{n}=a_{1} \alpha ^{n}+a_{2} \beta^{n}+a_{3} \gamma^{n} \quad for\, all\, n\geq 0 \\
\end{equation*}
The initial conditions $N_0=0,N_1=1$ and $ N_2=1 $ imply that
\begin{equation*}
  a_{1}=\textstyle\frac{\alpha}{(\alpha-\beta)(\alpha-\beta)},
  a_{2}=\textstyle\frac{\beta}{(\beta-\gamma)(\beta-\alpha)},
  a_{3}=\textstyle\frac{\gamma}{(\gamma-\alpha)(\gamma-\beta)}
\end{equation*}
The above  like Bient formula can also be written as
\begin{equation*}
  N_{n}=c_{\alpha}\alpha^{n+2}+c_{\beta}\beta^{n+2}+c_{\gamma}\gamma^{n+2}
\end{equation*}
where,
\begin{equation*}
    c_{t}=\textstyle\frac{1}{t^{3}+2} \quad , t\in \{\alpha,\beta,\gamma\}
\end{equation*}
It 's easy to verify the following inequalities approximations
\begin{equation*}
  1.45<\alpha<1.5
\end{equation*}
\begin{equation*}
  0.82<\lvert \gamma \vert =\lvert \beta \lvert <0.83
\end{equation*}
\begin{equation}
  5<c_{\alpha}^{-1}<5.15
  \label{equ:1}
\end{equation}
\begin{equation*}
  |c_{\beta}|\simeq0.4075
\end{equation*}
\begin{equation*}
  |\xi(n)|<\frac{1}{2} \quad where\, \xi(n)= c_{\beta}\beta^{n+2}+c_{\gamma}\gamma^{n+2}
\end{equation*}
By induction over $n$, it is easy to prove the relation between Narayana and $ \alpha $ \\
\begin{equation}
  \alpha^{n-2} \leq\ N_{n}\leq \alpha^{n-1}\quad for\, all\,  n \geq 0 \\
  \label{equ:3}
\end{equation}
We have
\begin{equation*}
    2^{l-1}\leq b^{l-1}\leq a \textstyle\frac{b^{l}-1}{b-1}=N_{n}N_{m}\leq \alpha^{n+m-2}\leq
    \alpha^{2n-2}\leq(1.5)^{2n-2}
\end{equation*}
\begin{equation*}
    l\leq(2n-2)\textstyle\frac{\log1.5}{\log2}+1 < 2n-1
\end{equation*}
and,
\begin{equation*}
 (1.45)^{n-2}<\alpha^{n-2}<N_{n}< N_{n}N_{m}=a \textstyle\frac{b^{l}-1}{b-1}<b^{l}<(50)^{l}
\end{equation*}
\begin{equation*}
    n<l\frac{\log10}{\log1.45}+2<11\,l+2
\end{equation*}
Similarly, we have
\begin{equation*}
2^{l_1-1}<b^{l_1-1} < \frac{b^{l_1}-1}{b-1}<a_1 a_2 \frac{(b^{l_1}-1)(b^{l_2}-1)}{(b-1)^2} =N_k<\alpha^{k-1}
\end{equation*}
\begin{equation*}
  l_{1}<(k-1) \frac{\log \alpha}{\log 2}+1<k
\end{equation*}
and
\begin{equation*}
  \alpha^{k-2}<N_k=a_1 a_2 \frac{(b^{l_1}-1) (b^{l_2}-1)}{(b-1)^2} <(b^{l_{2}}-1)^2<b^{2l_2}<50^{2l_2}
\end{equation*}
\begin{equation}
  \begin{split}
     k & < 2 l_2 \frac{\log 50}{\log \alpha}+2 \\
       &< 22 l_2+2
       \label{eq:a}
  \end{split}
\end{equation}
\subsection{Linear forms in logarithms of real algebraic number}
Let $\psi$ be an algebraic number of degree $ d$ with minimal polynomial over $\mathbb{Z}$
\begin{equation*}
 f(X)=a_{0}\prod_{i=1}^{d}(X-\psi^{(i)}).
\end{equation*}
where $a_{0}>0$ is leading coefficient, and  $\psi^{(i)}$'s are the conjugates of $\psi$. The logarithmic height of
$\psi$  $[\cite{Di_An},Def. \,2.2.8]$ is defined by
\begin{equation*}
  h(\psi)=\frac{1}{d} (\log a_{0}+\sum_{i=1}^{d}\log\max\{|\psi^{(i)}|,1\}).
\end{equation*}
 and the following properties hold:
\begin{equation}
  \begin{split}
     h(\psi \pm \gamma) & \leq h(\psi)+h(\gamma)+\log 2 \\
     h(\psi \gamma^{\pm1} ) & \leq h(\psi)+h(\gamma) \\
       h(\psi^{s})& =|s| h(\psi)  \quad \quad (s\in \mathbb{Z})
  \end{split}
  \label{prop:1}
\end{equation}

\begin{theorem}
[{\color{blue}(Matveev)},\cite{Dujella}] Let $\psi_{1},\dots \psi_{t}$ be positive real algebraic numbers, $\mathbb{K}$ be a
number field of degree D over over $\mathbb{Q}$, and $r_1,\dots ,r_t $ integers. Let
\begin{equation*}
  \Lambda=\psi^{r_1}_{1}\cdots\psi^{r_t}_{t}
\end{equation*}
let $B \geq \max\{|r_1|,\cdots |r_t|\}$ and $A_{j}\geq \max \{Dh(\psi_{j}),|\log \psi_{j}|,0.16 |\}$ if $\Lambda\neq
0$, then $\log|\Lambda|>-1.4 \times 30^{t+3} \times t^{4.5} \times D^{2} (1+\log D )(1+\log B) A_{1}\cdots A_{t} $.
\label{th:bd1}
\end{theorem}
\begin{lemma}
  $[\cite{Ssnchez},Lemma \,7]$ If $m\geq1$, $T>(4 m^{2})^{m}$ and $ T>\frac{x}{\log^{m} x}$, then $x<2^{m} T
  \log^{m} T$.
  \label{l:1}
\end{lemma}
 This $\mathbf{lemma}$  will be used to reduce the upper bound for variables, and we will define
 $\|X\|=\min\{|X-n|:n\in  \mathbb{Z}\}$ be the dinstance from $X $ to the nearest integer.
\begin{lemma}
({\color{blue}(Dujella- peth\"{o})} ,$\cite{Di_An},Lemma \, 2.3.1)$
 Let M be a positive integer such that $q>6M$, since $\frac{p}{q}$ is a convergent of the irrational number $\tau$,
 let A,B, and $\mu$  be some real numbers with $A>0$,$  B>1$ and $\epsilon=\|\mu q\|-M \|\tau q\|$. if $\epsilon>0$,
 then there is no solution to the inequality
  \begin{equation*}
    0<|u\tau -v+\mu|<A B^{-w}
  \end{equation*}
  in positive integers u ,v and w with
  \begin{equation*}
    u\leq M \quad  and \quad w\geq\frac{\log(Aq/\epsilon)}{\log B}
  \end{equation*}
  \label{lem:2}
\end{lemma}

\begin{lemma}
 (i)  [{\color{blue}(Legender)} ,$\cite{Di_An},Theorem\, 1.3.3]$ Let $\tau$ be an irrational number such that
  \begin{equation*}
    |\tau -\frac{x}{y}|<\frac{1}{2 y^{2}}
  \end{equation*}
  then $\frac{x}{y}$ is a convergent of $\tau$.\\
  \\
  (ii) If $y<q_{k+1}$ then
  \begin{equation*}
\frac{1}{(g+2) y^{2}} < |\tau -\frac{x}{y}|
  \end{equation*}
   $g=max\{g_{i}\, : j\leq k+1\}$.
   \label{le:3}
\end{lemma}

\section{Proof of theorem \ref{th:u}}
\subsection{Bounding on $l_{1}$}
From equation{\color{blue}(\ref{equ:r})}, we obtain that
\begin{equation*}
  c_{\alpha} \alpha^{k+2}-\frac{a_1 a_2 b^{l_1+l_2}}{(b-1)^2}=-\xi(k)-\frac{a_1 a_2 b^{l_1}}{(b-1)^2}-\frac{a_1 a_2
  b^{l_2}}{(b-1)^2}+\frac{a_1 a_2}{(b-1)^2}
\end{equation*}
Taking absolute values in the above equation, using inequalities {\color{blue}(\ref{equ:1}),(\ref{equ:3})} and
dividing both sides by $|\frac{a_1 a_2 b^{l_1+l_2}}{(b-1)^2}|$, we get
\begin{equation*}
  \begin{split}
     | c_{\alpha} \alpha^{k+2}-\frac{a_1 a_2 b^{l_1+l_2}}{(b-1)^2}| & < \frac{1}{2}+b^{l_1}+b^{l_2}+1 \\
       & <\frac{3}{2} +2 b^{l_2}\\
      |\frac{c_{\alpha} \alpha^{n+2} (b-1)^{2}}{a_1 a_2 b^{l_1+l_2}}-1| & <\frac{3 (b-1)^{2}}{2 a_1 a_2
      b^{l_1+l_2}}+\frac{2 (b-1)^2}{a_1 a_2 b^{l_1}} \\
       & <\frac{3 (b-1)^{2}}{b^{l_1}}+\frac{2(b-1)^2}{b^{l_1}} \\
       & <\frac{3 b^{2}}{b^{l_1}}+\frac{2 b^2}{b^{l_1}} \\
       & <\frac{5}{b^{l_1-2}}
  \end{split}
\end{equation*}
Put
\begin{equation*}
  \Lambda_{3}=\frac{c_{\alpha} \alpha^{n+2} (b-1)^{2}}{a_1 a_2 b^{l_1+l_2}}-1
\end{equation*}
we have
\begin{equation}
  |\Lambda_{3}|<\frac{5}{b^{l_1-2}}
  \label{eq:gt}
\end{equation}
and $\log|\Lambda_{3}| <\log 5 -(l_1-2)\log b$ Now, we apply matveev theorem, where
\begin{equation*}
  \begin{array}{ccc}
    \psi_1=\alpha & \psi_2=b & \psi_3=\frac{c_{\alpha} (b-1)^2}{a_1 a_2} \\
r_1=(k+2)  & r_2=-(l_1+l_2) & r_3=1
  \end{array}
\end{equation*}
Similarly we can prove that $ \Lambda_{3}\neq 0$, moreover using   properties of logarithmic height (\ref{prop:1}),
we obtain
\begin{equation*}
  \begin{split}
     h(\psi_3) & < h(c_{\alpha})+h(\frac{b-1}{a_1})+ h(\frac{b-1}{a_2})\\
       &<\frac{\log 31}{3}+2 \log (b-1) \\
       & < 3\log b
  \end{split}
\end{equation*}
Thus, we can take $A_1=\log \alpha$ ,$A_2=3 \log b$ , $A_3=9 \log b$ , $B=22 l_2+4$ since $k<22 l_2+2$  and
$\mathbb{K}=\mathbb{Q}(\alpha)$ thus $D=3$, and then from theorem (\ref{th:bd1}) we get
\begin{equation*}
  \log \Lambda_{3}>-1.4 \cdot 30^{6} \cdot 3^{4.5} \cdot 3^{5} (1+\log 3) (1+\log(22 l_2+4)) \log \alpha \log^2 b
\end{equation*}
Now we compare the lower bound for $\log \Lambda_{3}$ with the upper bound of $ \log \Lambda_{3}$. Since $(1+\log(22
l_2+4) )<8 \log (l_2) $ for all $ l_2 \geq 2$, a computer
search with Mathematica gives us that
\begin{equation}
  l_1< 3 \times 10^{14} \log l_2 \log b
  \label{eq:r}
\end{equation}
\subsection{Bounding on $l_{2}$}
Let
\begin{equation*}
  \begin{split}
     \frac{N_k}{\dfrac{a_1 (b^{l_1}-1)}{b-1}}& =\frac{a_2 (b^{l_2}-1)}{b-1} \\
     \frac{c_{\alpha} \alpha^{k+2} (b-1)}{a_1 (b^{l_1}-1)} -\frac{a_2 b^{l_2}}{b-1} &=\frac{-\xi(k) (b-1)}{a_1
     (b^{l_1}-1)}-\frac{a_2}{b-1}
  \end{split}
\end{equation*}
Taking absolute values in the above equation and
dividing both sides by $|\frac{a_2 b^{l_2}}{b-1}|$, we get
\begin{equation*}
  \begin{split}
     |\frac{c_{\alpha} \alpha^{k+2} (b-1)}{a_1 (b^{l_1}-1)} -\frac{a_2 b^{l_2}}{b-1}| & <\frac{(b-1)}{2 a_1
     (b^{l_1}-1)} +1\\
      |\frac{\alpha^{k+2} c_{\alpha} b^{-l_2} (b-1)^2}{a_1 a_2 (b^{l_1}-1)}-1| & <\frac{(b-1)^2}{a_1 a_2 b^{l_2}
      (b^{l_1}-1)}+\frac{b-1}{a_2 b^{l_2}} \\
       & <\frac{(b-1)^2}{b^{l_2}}+\frac{b-1}{b^{l_2}} \\
       & < \frac{b^2}{b^{l_2}}+\frac{b}{b^{l_2}}
  \end{split}
\end{equation*}
\begin{equation}
  |\frac{\alpha^{k+2} c_{\alpha} b^{-l_2} (b-1)^2}{a_1 a_2 (b^{l_1}-1)}-1)|<\frac{2}{b^{l_2-2}}
  \label{eq:m}
\end{equation}
Put $ \Lambda_4=\frac{\alpha^{k+2} c_{\alpha} b^{-l_2} (b-1)^2}{a_1 a_2 (b^{l_1}-1)}$, we have
\begin{equation}\label{eq:i}
  \log|\Lambda_4|<\log 2-(l_2-2) \log b
\end{equation}
Now, we apply matveev theorem ${\color{blue}(\ref{th:bd1})}$, where
\begin{equation*}
  \begin{array}{ccc}
    \psi_1=\alpha &\psi_2= b & \psi_3=\frac{c_{\alpha} (b-1)^2}{a_1 a_2 (b^{l_1}-1)} \\
    r_1=k+2 & r_2=-l_2 & r_3=1
  \end{array}
\end{equation*}
Similarly we can prove that $|\Lambda_4|\neq 0$, moreover using properties of logarithmic
height (\ref{prop:1})
\begin{equation*}
  \begin{split}
     h(\psi_3) & < h(c_{\alpha})+h(\frac{b-1}{a_1})+h(\frac{b-1}{a_2})+h(b^{l_1}-1) \\
       & < \frac{\log 31}{3}+2 \log(b-1)+l_1 \log b\\
       & < 3 \log b+l_1 \log b
  \end{split}
\end{equation*}
thus, we can take $A_1=\log \alpha$ ,$A_2=3 \log b$ , $A_3=3(4 \log b +l_1 \log b)$ and $B=22 l_2+4$
\begin{equation}
  \log \Lambda_{4}>-1.4 \cdot 30^{6} \cdot 3^{4.5} \cdot 3^{4} \log \alpha (1+\log 3) (1+\log(12 l_2+2))  (4 \log
  b+l_1 \log b)
  \label{eq:ii}
\end{equation}
from ${\color{blue}(\ref{eq:r})}$  ,${\color{blue}(\ref{eq:i})}$ and $(\ref{eq:ii})$ we deduce that
\begin{equation*}
   l_2  <2 \times 10^{28} \log b \log^2 l_2
\end{equation*}
Now we apply lemma $(\ref{l:1})$, since $ 2 \times 10^{28} \log^{2}(l_2)  \log b> (16)^2 $, we obtain
\begin{equation*}
  \begin{split}
    \frac{l_2}{\log^2 l_2}   & < 2 \times 10^{28} \log b \\
      l_2 & < 2^2 \cdot 2 \cdot 10^{28}\log b (\log(2 \times 10^{28} \log b))^2 \\
       & <10^{29} \log(b) (66+\log \log b)^2 \\
       &< 10^{33} \log^3 b
  \end{split}
\end{equation*}
since $(66+\log\log b )^{2}< 95^2 \log^{2} b  $ for every $ b\geq 2$. from {\color{blue}(\ref{eq:a})}, we find that
$k<2.3 \times 10^{34} \log^3 b$.

\subsection{Reduction of  The upper bound on $l_1$}
Let $z_3=(n+2) \log \alpha -(l_1+l_2) \log b +\log \frac{(b-1)^2 c_{\alpha}}{a_1 a_2}$ ,\ if $z_3>0$  then
$z_3<|e^{z_3} -1|$ and $ |z_3|<2 |e^{z_3} -1| \, if z_3<0$,Thus in both side we have, $ |z_3|<2 |e^{z_3} -1|$. By
substituting into the equation ${\color{blue}(\ref{eq:gt})}$, dividing both by $\log b$, we have
\begin{equation*}
  \begin{split}
  |(k+2) \log \alpha-(l_1+l_2)\log b+\log (\frac {(b-1)^2 c_{\alpha}}{a_1 a_2})| & <\frac{10}{b^{l_1-2}}\\
  |(k+2)\frac{\log \alpha}{\log b}-(l_1+l_2)+ \frac{\log(\dfrac{(b-1)^2 c_{\alpha}}{a_1 a_2})}{\log b}|  &
  <\frac{10}{\log(b) b^{l_1-2}}
  \end{split}
  \end{equation*}
  \begin{equation}
  |(k+2)\frac{\log \alpha}{\log b}-(l_1+l_2)+ \frac{\log(\dfrac{(b-1)^2 c_{\alpha}}{a_1 a_2})}{\log b}|
  <\frac{15}{b^{l_1-2}}
  \label{eq:f}
\end{equation}
Since $\frac{1}{\log 2}=1.4427$. Let $\tau =\frac{\log \alpha}{\log b} , \mu=\frac{\log(\dfrac{(b-1)^2
c_{\alpha}}{a_1 a_2})}{\log b}$ and $M=1.3 \times 10^{34}\log^{3}b$, at all $  b\in \{2,3,\cdots,50\}$ and $ a_1
,a_2\in \{1,\cdots,b-1\}$, a computer search with Mathematica find that $\varepsilon >0$ for all, so we apply lemma
{\color {blue}(\ref{lem:2})}, let $A=15$ and $B=b$, we can say that if the inequality {\color {blue}(\ref{eq:f})}
has a solution then $l_1-2\leq \max(\frac{\log(\dfrac{A q_k}{\varepsilon})}{\log B})\leq120$, hence $l_1\leq122$.
\subsection{Reduction of The upper bound on $l_2$}
Let $z_4=(k+2) \log \alpha -l_2 \log b +\log \frac{c_{\alpha} (b_-1)^2}{a_1 a_2 (b^{l_1}-1)}$ ,\ if $z_4>0$  then
$z_4<|e^{z_4} -1|$ and $ |z_4|<2 |e^{z_4} -1| \, if z_4<0$, thus in both side we have, $ |z_4|<2 |e^{z_4} -1|$. By
substituting into the equation ${\color{blue}(\ref{eq:m})}$ and dividing both by $\log b$, we have
\begin{equation*}
  \begin{split}
     |(k+2) \frac{\log \alpha}{\log b}-l_2+\frac{\log(\dfrac{c_{\alpha} (b-1)^2}{a_1 a_2 (b^{l_1}-1)})}{\log b}| & <
     \frac{4}{\log b \,b^{l_2-2}} \\
       & <\frac{6}{ b^{l_2-2}}
  \end{split}
\end{equation*}
Let $\tau =\frac{\log \alpha}{\log b} , \mu=\frac{\log(\dfrac{c_{\alpha} (b-1)^2}{a_1 a_2 (b^{l_1}-1)})}{\log b}$
and $M=1.3 \times 10^{34}\log^{3}b_1$  at all $  b\in \{2,3,\cdots,10\}$, $ a_1,a_2\in \{1,\cdots,b-1\}$,and $l_1
\in \{1,\cdots,122\}$, a computer search with Mathematica founds that $\varepsilon >0$ for all, so we apply lemma
{\color{blue}(\ref{lem:2})}, let $A=6$ and $B=b$, we can say that if the inequality {\color{blue} (\ref{eq:f})} has
a solution then $l_2-2\leq \max(\frac{\log(\dfrac{A q_k}{\varepsilon})}{\log B})\leq131$ , hence $l_2\leq133$, then
$k<1598$.
%------

\section{Proof of theorem  \ref{th:1}}
\subsection{Bounding on  m}
From equation (\ref{equ:2}), we obtain that
\begin{equation*}
\end{equation*}
\begin{equation*}
  c^{2}_{\alpha}\alpha^{n+m+4}-\frac{a b^{l}}{b-1}=-\xi(m)c_{\alpha}\alpha^{n+2}
  -\xi(n)c_{\alpha}\alpha^{m+2}-\xi(n)\xi(m)-\frac{a}{b-1}
\end{equation*}
Taking absolute values in the above equation, using inequalities {\color{blue}(\ref{equ:1})} and dividing both sides
by $|c^{2}_{\alpha}\alpha^{n+m+4}|$, one gets
\begin{equation*}
  \begin{split}
  \Big |c^{2}_{\alpha}\alpha^{n+m+4}-\frac{a b^{l}}{b-1} \Big|& <
  \frac{c_{\alpha}\alpha^{n+2}}{2}+\frac{c_{\alpha}\alpha^{m+2}}{2}+\frac{5}{4} \\
    \Big |1-\frac{a b^{l}}{c^{2}_{\alpha}\alpha^{n+m+4}(b-1)} \Big | &<
    \frac{1}{2c_{\alpha}\alpha^{m+2}}+\frac{1}{2c_{\alpha}\alpha^{n+2}}+\frac{5}{4 c^{2}_{\alpha}\alpha^{n+m+4}} \\
       & <\frac{1}{c_{\alpha}\alpha^{m+2}}+\frac{5}{4 c^{2}_{\alpha}\alpha^{m+2}}  \\
       & < \frac{39}{\alpha^{m}}
  \end{split}
\end{equation*}
Put
\begin{equation*}
  \Lambda_{1}:=\frac{a b^{l}}{c^{2}_{\alpha}\alpha^{n+m+4}(b-1)}-1
\end{equation*}
we have
\begin{equation}
  |\Lambda_{1}|< \frac{39}{\alpha^{m}} \quad and \,  \log|\Lambda_{1}|<\log(39)-m\log(\alpha)
  \label{equ:4}
\end{equation}
Now, we apply the Matveev theorem, where
\begin{equation*}
  \begin{array}{ccc}
    \psi_{1}=\alpha & \psi_{2}=b & \psi_{3}=\frac{a}{c^{2}_{\alpha} (b-1)}\\
     r_1=-(n+m+4) & r_2=l  & r_3=1
  \end{array}
\end{equation*}
First, we show that $\Lambda_{1} \neq 0$. If $\Lambda_{1}= 0$, then $ \frac{a
b^{l}}{b-1}=c^{2}_{\alpha}\alpha^{n+m+4}$. Consider the automorphism $\sigma(c_{\alpha})=c_{\beta}$.Then
$|c^{2}_{\beta}\beta^{n+m+4}|<|c^{2}_{\beta}|<1$, while the right-hand side is greater than 4
which is a contradiction, moreover using properties of logarithmic height (\ref{prop:1}),
we obtain
\begin{equation*}
  h(\psi_{1})=\frac{\log(\alpha)}{3},  h(\psi_{2})=\log(b)
\end{equation*}
\begin{equation*}
  \begin{split}
 h(\psi_{3}) & <h(\frac{a}{b-1})+h(c^{2}_{\alpha}) \\
              &< \log(b-1)+\frac{2 \log(31)}{3} \\
              &<\log(b)+3.4\log(b)\\
              &< 4.5\log(b)
\end{split}
\end{equation*}
 since the minimal polynomial of $c_{\alpha}$ is given by $31 x^{3}-31 x^{2}+10x-1$. We take
$B=2n+4$, $A_{1}=\log(\alpha)$, $ A_{2}=3 \log(b)$, $ A_{3}=13.5\log(b)$, we take
$\mathbb{K}=\mathbb{Q}(\alpha),$ thus $D=3$.\\
Now from  theorem (\ref{th:bd1}), we get the following
\begin{equation*}
  \log|\Lambda_{1}|>-1.4\cdot 30^{6}\cdot 3^{4.5} \cdot 3^{3}\cdot 13.5 \, (1+\log(3))(1+\log(2n+4)) \, \log(\alpha)
  \log^{2}(b)
\end{equation*}
Now we compare the lower bound for  $\log|\Lambda_{1}|$ with the upper bound of $ \log|\Lambda_{1}|$.
Since $1+log(2n+4)<5log(n)$ for all $n\geq3$, a computer search with Mathematica gives us that
\begin{equation}
  m<1.7\times 10^{15} \log(n) \log^{2}(b)
  \label{equ:p}
\end{equation}
\subsection{Bounding on n}
Let
\begin{equation*}
   \begin{split}
      N_{n} & =\frac{a}{N_{m}} \frac{b^{l}-1}{b-1} \\
      c_{\alpha} \alpha^{n+2}-\frac{a b^{l}}{N_{m} (b-1)}  & =-\xi(n)-\frac{a}{N_{m} (b-1)}
   \end{split}
\end{equation*}
Taking absolute values in the above equation, using inequalities {\color{blue}(\ref{equ:1}),(\ref{equ:3}) }and
dividing both sides by $|c_{\alpha}\alpha^{n+2}|$, we get
\begin{equation}
  \begin{split}
    \Big |c_{\alpha}\alpha^{n+2}-\frac{a b^{l}}{N_{m} (b-1)} \Big| &<|\xi(n)|+|\frac{a}{N_{m} (b-1)} |\\
       &< \frac{1}{2}+\frac{1}{\alpha^{m-2}} \\
 \Big |  1-\frac{a b^{l}}{ N_{m} c_{\alpha}\alpha^{n+2}(b-1)} \Big|    & <\frac{1}{2
 c_{\alpha}\alpha^{n+2}}+\frac{1}{c_{\alpha} \alpha^{n+m}} \\
 &<\frac{1}{2c_{\alpha} \alpha^{n}}+\frac{1}{c_{\alpha} \alpha^{n}}\\
       & < \frac{11}{\alpha^{n}}
  \end{split}
\end{equation}
Put
\begin{equation*}
  \Lambda_{2}:=\frac{a b^{l}}{ N_{m} c_{\alpha}\alpha^{n+2}(b-1)}-1
\end{equation*}
we have
\begin{equation}
  |\Lambda_{2}|<\frac{11}{\alpha^{n}}
  \label{equ:e}
\end{equation}
and $ \log|\Lambda_{2}|<\log(11)-n\log(\alpha)$.
Now, we apply matveev theorem (\ref{th:bd1}), where
\begin{equation*}
  \begin{array}{ccc}
    \psi_{1}=\alpha & \psi_{1}=b & \psi_{1}=\frac{a}{N_{m} c_{\alpha}(b-1)} \\
    r_1=-(n+2) & r_2=l &r_3= 1
  \end{array}
\end{equation*}
Similarly we can prove that $\Lambda_{2}\neq 0$, moreover using   properties of logarithmic height (\ref{prop:1}),
we obtain
\begin{equation*}
\begin{split}
 h(\psi_{3}) & <h(\frac{a}{b-1})+h(c_{\alpha})+h(N_{m}) \\
              &< \log(b-1)+\frac{\log(31)}{3}+m\log(\alpha)\\
              &<\log(b)+1.2\log(b)+m\log(\alpha)\\
              &<2.3\log(b)+m\log(\alpha)
\end{split}
\end{equation*}
 we take $B=2n+2$, $ A_{1}=\log(\alpha),\, A_{2}=3 \log(b)\, , A_{3}=3(2.3
 \log(b)+m\log(\alpha))$, $\mathbb{K}=\mathbb{Q}(\alpha)$ thus $D=3$, from theorem (\ref{th:bd1}) we get
\begin{equation*}
\log |\Lambda_{2}|> -1.4 \cdot 30^{6} \cdot 3^{4.5} \cdot 3^{4}\, \log(\alpha)\, \log(b)\, (1+\log(3))\,
(1+\log(2n+2))\, (2.3\,\log(b)+m\,\log(\alpha)).
\end{equation*}
Now we compare the lower bound for  $\log|\Lambda_{2}|$ with the upper bound of $ \log|\Lambda_{2}|$
 and using {\color{blue}(\ref{equ:p})}, a computer search with Mathematica gives us that
\begin{equation*}
  \begin{split}
      n & <7.6\times 10^{28}\, \log^{2}n \, \log^{3}b \\
     \frac{n}{\log^{2}(n)}  &  < 7.6\times 10^{28} \log^{3}b \\
  \end{split}
\end{equation*}
Now we apply lemma {\color{blue}(\ref{l:1}) }, since $7.6\times 10^{28} \log^{3}(b)>(16)^{2}$, we obtain
\begin{equation}
  \begin{split}
                        n   & <2^{2} \cdot 7.6 \cdot 10^{28} \log^{3}(b) (\log(7.6 \times 10^{28} \log^{3}b))^{2} \\
                           & < 3.04\times 10^{29} \log^{3}b (66.6+3 \log \log b )^2 \\
                           &< 3.04\times 10^{29} \log^{3}b (96.1\,\log b+3 \log b )^2\\
                           &<2.99 \times10^{33} \log^{5}b
\end{split}
\end{equation}
since \,$\log \log b< \log b $ \, for every $ b\geq 2$ and $\frac{1}{\log 2}\simeq 1.4427$.
\subsection{ Reduction of The upper bound on m}
Let $z_1=l \log(b)-(n+m+4) \log\alpha +\log(\frac{a}{(b-1) c^{2}_{\alpha}})$ ,\ if $z_1>0$  then $z_1<|e^{z_1} -1|$
and $ |z_1|<2 |e^{z_1} -1| \, if\, z_1<0$, thus in both side we have, $ |z_1|<2 |e^{z_1} -1|$. By substituting into
the equation (\ref{equ:4}), we have
\begin{equation*}
  |l \log b-(n+m+4)\log(\alpha)+\log(\frac{a}{(b-1) c^{2}_{\alpha}}) |<\frac{78}{\alpha^{m}}
\end{equation*}
Dividing this inequality by $|\log \alpha|$, we get
\begin{equation}
  |l \frac{\log b}{\log \alpha} -(n+m+4)+\frac{\log(\dfrac{a}{c^{2}_{\alpha} (b-1)})}{\log
  \alpha}|<\frac{210}{\alpha^{m}}
  \label{equ:8}
\end{equation}
Let $ \tau =\frac{\log(b)}{\log \alpha} , \mu =\frac{\log(\dfrac{a}{c^{2}_{\alpha} (b-1)})}{\log \alpha}$ and
$M=5.98\times10^{33} \log^{5} b$. For all $  b\in \{2,3,\cdots,50\}$ and $ a\in \{1,2,\cdots,b-1\}$, we need to
calculate a convergent $ \frac{p_{k}}{q_{k}}$ such that $q_{k}>6 M$, furthermore computing
$\varepsilon=\|\mu q_{k}\|-M\|\tau q_{k}\| $, a computer search with Mathematica find that $\varepsilon >0$ for all,
so we can apply lemma {\color{blue}(\ref{lem:2})}, let $A=210$, and $B=\alpha$, we can say that if the inequality
{\color{blue}(\ref{equ:8})} has a solution then $m\leq \max \Bigg ( \frac{\log(\dfrac{A q_{k}}{\varepsilon})}{\log
B} \Bigg ) \leq 261$.
\subsection{Reduction of The upper bound on n}
Let $z_{2}=l\log b-(n+2)\log \alpha +\log(\frac{a}{N_{m} c_{\alpha} (b-1)})$ , substituting into the equation
{\color{blue} (\ref{equ:e})}, we have
\begin{equation}
 \Big |l \frac{\log b}{\log \alpha}-(n+2)+\frac{\log(\dfrac{a}{N_{m} c_{\alpha} (b-1)})}{\log
 \alpha}\Big|<\frac{32}{\alpha^{n}}
 \label{equ:10}
\end{equation}
Let $\tau=\frac{\log b}{\log\alpha}$ , $ \mu=\frac{\log(\dfrac{a}{N_{m} c_{\alpha} (b-1)})}{\log \alpha} $ and
$M=5.98\times10^{33} \log^{5} b $ , at all $  b\in \{2,3,\cdots,50\}$ , $ a\in \{1,2,\cdots,b-1\}$ and $m\in
\{3,\cdots,261\}$, a computer search with Mathematica find that $\varepsilon >0$ for all except $(b,a,m)=$ $\{
(b,b-1,3) for\, all\, b={2,\cdots,50} \}$, in addition
to$\{(2,1,4),(2,1,6),(3,2,5),(3,2,8),(4,3,6),(6,5,7),(9,8,8),(13,12,9)\\,(19,18,10),(28,27,11),(41,40,12)\}$. We apply lemma
{\color{blue}(\ref{lem:2})} in case $\varepsilon >0$, let $A=32$ and $B=\alpha$, we can say that if the inequality
{\color{blue}(\ref{equ:10})} has a solution then $n\leq \max (\frac{\log(\dfrac{A q_k}{\varepsilon})}{\log
B})\leq290$, in other cases we apply Lemma {\color{blue}(\ref{le:3})},
\begin{equation}
 \Big | \frac{\log b}{\log \alpha}-\frac{(n+2)-\frac{\log(\dfrac{a}{N_{m} c_{\alpha} (b-1)})}{\log
 \alpha}}{l}\Big|<\frac{32}{ \alpha^{n} l}
 \label{eq:100}
\end{equation}
 now assume that $ n $ is so large the right hand side of the inequality {\color{blue}(\ref{eq:100})} is smaller
 than $\frac{1}{2 l^{2}}$ holds if $\alpha^{n}>64 l$, which by Lemma {\color{blue}(\ref{le:3})}, implies that the
 fraction $\frac{\log b}{\log \alpha}$ is a convergent of $\frac{(n+2)-\frac{\log(\dfrac{1}{N_{m} c_{\alpha}
 })}{\log \alpha}}{l}$, since in all case $a=b-1$, for each $( b,a ,m)$ which have $\varepsilon <0$, we calculate
 the continued fraction expantion of $\tau $ and find \\$g=max\{g_{i} \, : j\leq k+1\}$. since
\begin{equation*}
\frac{1}{(g+2) l^2}<\big|\frac{\log b}{\log
\alpha}-\frac{(n+2)-\dfrac{\log(\dfrac{a}{N_{m}c_{\alpha}(b-1)})}{\log\alpha}}{l} \big|<\frac{32}{\alpha^{n} l}
\end{equation*}
\begin{equation*}
  \begin{split}
     \alpha^{n} & < 32(g+2) l \\
     n  & <\frac{\log(32(g+2) l)}{\log \alpha} \\
       & <\frac{\log(32\times 5.98\times 10^{33} \log^{5} b (g+2))}{\log \alpha}
  \end{split}
\end{equation*}
we found $n\leq239$, therefore $ n\leq290$ in both cases.\\
\\
We conclude all solutions $(n,m,l,a,b)$ to the Diophantine equation (\ref{equ:2}) $3\leq m\leq n,2\leq b\leq
50,1\leq a\leq b-1 $ and $ l \geq 2$, reduce to the rang $3\leq n\leq 264$, with the help of Mathematica, we compute
all solution in specified range, we conclude theorem (\ref{th:1}).

%------
\addcontentsline{toc}{chapter}{Bibliography}
\markboth{\textbf{\textit{Bibliography}}}{
\end{document}